\documentclass[11pt,reqno]{amsart}
\usepackage[T1]{fontenc}
\usepackage{lmodern}
\usepackage{amssymb}
\usepackage[hidelinks]{hyperref}

\newtheorem{theorem}{Theorem}
\newtheorem{proposition}[theorem]{Proposition}
\newcommand{\Aut}{\operatorname{Aut}}
\newcommand{\Tr}{\operatorname{Tr}}

\title[Centerless groups and rack automorphisms]{A centerless group that is not the automorphism group of a rack}
\author{Amir Saki}
\address{Universit\'e du Qu\'ebec \`a Trois-Rivi\`eres}
\email{amir.saki.math@gmail.com}
\date{}

\begin{document}

\begin{abstract}
We construct a centerless group of order \(1280\) that is not isomorphic to
the automorphism group of any rack. This answers a question of Militaru
for both racks and quandles. The obstruction comes from an outer
automorphism that preserves the simultaneous conjugacy class of each
pair consisting of a point stabilizer and a right translation in any
possible rack action.
\end{abstract}

\maketitle

\section{Introduction}

Let \(X\) be a rack with operation \(*\). For each \(x\in X\), let \(R_x\)
denote the right translation defined by \(R_x(y)=y*x\). Thus each \(R_x\)
is a bijection, and right self-distributivity implies that \(R_x\) is an
automorphism of \(X\). A rack is a quandle if \(x*x=x\) for each \(x\in X\).

Militaru proved that some finite groups, including the quaternion group
of order eight, are not isomorphic to the automorphism group of any
rack. He asked whether every finite group with trivial center is the
automorphism group of a finite rack or of a finite quandle
\cite[Question~2]{Militaru}.

\begin{theorem}\label{thm:main}
There exists a centerless group of order \(1280\) that is not isomorphic
to the automorphism group of any rack.
\end{theorem}

\section{An obstruction from point stabilizers}

For a group \(G\) and \(s\in G\), write \(C_G(s)\) for the centralizer of
\(s\). For \(H\le G\), write \(N_G(H)\) for the normalizer of \(H\).
We say that two pairs \((H,s)\) and \((H',s')\) are conjugate in \(G\) if
there exists \(a\in G\) with \(H'=aHa^{-1}\) and \(s'=asa^{-1}\).

\begin{proposition}\label{prop:obstruction}
Let \(G\) be a centerless group. Assume that \(G\) has an outer
automorphism \(\alpha\) such that \((\alpha(H),\alpha(s))\) is conjugate
to \((H,s)\) for every \(s\in G\) and every \(H\le C_G(s)\) with
\(N_{C_G(s)}(H)=H\). Then \(G\) is not isomorphic to the automorphism group
of any rack.
\end{proposition}

\newpage
\begin{proof}
Assume that \(G=\Aut(X)\) for a rack \(X\), and let \(G\) act on \(X\)
by evaluation. For \(x\in X\), let \(G_x\) denote its stabilizer.
Militaru's lemma \cite[Lemma~2.1]{Militaru} gives
\[
  G_x\le C_G(R_x)
  \quad\text{and}\quad
  N_G(G_x)\cap C_G(R_x)\le G_x Z(G).
\]
Since \(Z(G)=1\), we have \(N_{C_G(R_x)}(G_x)=G_x\).

For each \(G\)-orbit in \(X\), choose a representative \(x\).
By hypothesis there exists \(a_x\in G\) such that
\[
  \alpha(G_x)=a_xG_xa_x^{-1}
  \quad\text{and}\quad
  \alpha(R_x)=a_xR_xa_x^{-1}.
\]
Define \(f\) on this orbit by \(f(gx)=\alpha(g)a_xx\) for each \(g\in G\).
The equality of stabilizers makes \(f\) a well-defined bijection on the
orbit. Applying the definition on each orbit gives
\(f(gy)=\alpha(g)f(y)\) for all \(g\in G\) and \(y\in X\).
For \(y=gx\), the identity \(R_{gx}=gR_xg^{-1}\) gives
\[
  R_{f(y)}
  =\alpha(g)a_xR_xa_x^{-1}\alpha(g)^{-1}
  =\alpha(R_y).
\]
It follows that \(f(R_y(z))=R_{f(y)}(f(z))\) for all \(y,z\in X\).
Thus \(f\in G\). The equality \(f(gy)=\alpha(g)f(y)\) now gives
\(\alpha(g)=fgf^{-1}\) for each \(g\in G\), a contradiction.
\end{proof}

\section{The group of order \(1280\)}

Let \(F=\mathbb F_{16}\). On \(P=F\times F\), define a multiplication by
\[
  (x,z)(y,w)=(x+y,\ z+w+x^2y^4).
\]
Since \((x,y)\mapsto x^2y^4\) is additive in each variable, the
multiplication is associative. Its identity is \((0,0)\), and
\((x,z)^{-1}=(x,z+x^6)\). Put \(Z=\{0\}\times F\). Then \(Z\) is central in
\(P\) and \(P/Z\) is the additive group of \(F\). With the convention
\([u,v]=uvu^{-1}v^{-1}\), we have
\begin{equation}\label{eq:commutator}
  [(x,z),(y,w)]=(0,x^2y^4+y^2x^4).
\end{equation}

Choose \(a\in F^\times\) of order five and let \(D=\langle d\rangle\)
be cyclic of order five. Since \(a^6=a\), the map
\((x,z)\mapsto(ax,az)\) is an automorphism of \(P\).
Let \(G=P\rtimes D\), where \(d\) acts by this map. Then
\(|G|=16^2\cdot5=1280\). For \(1\le i\le4\), the action of \(d^i\) fixes
only \((0,0)\) in \(P\), so \(C_P(d^i)=1\). In particular \(C_G(d)=D\).
No nonidentity element of \(D\) centralizes \(P\), and hence \(Z(G)=1\).

Define an automorphism \(\alpha\) of \(G\) by
\[
  \alpha(d)=d,\qquad \alpha(x,z)=(x,z+x)\quad ((x,z)\in P).
\]
The map on \(P\) preserves the multiplication and commutes with the
action of \(D\). It has order two and induces the identity on \(P/Z\).
If \(\alpha\) were conjugation by \(g\in G\), then \(g\) would belong to
\(C_G(d)=D\). The action of \(D\) on \(P/Z\) is faithful, so \(g=1\), a
contradiction. Thus \(\alpha\) is outer.

\newpage
\begin{proof}[Proof of Theorem~\ref{thm:main}]
We check the hypothesis of Proposition~\ref{prop:obstruction} for \(G\)
and \(\alpha\). Let \(s\in G\), and let \(H\le C_G(s)\) satisfy
\(N_{C_G(s)}(H)=H\).

First assume that \(s\ne1\). For \(1\le i\le4\), the \(P\)-conjugacy class
of \(d^i\) has size \(|P|\) and is contained in \(Pd^i\). Consequently
each element of \(G\setminus P\) is conjugate to a nonidentity element
of \(D\). If \(s\in G\setminus P\), then \(C_G(s)\) is conjugate to \(D\).
If \(s\in P\setminus\{1\}\), then \(C_G(s)\le P\), since any element
of \(G\setminus P\) has trivial centralizer in \(P\).
Thus \(C_G(s)\) is a finite \(p\)-group, for \(p=2\) or \(p=5\).
Every proper subgroup of a finite \(p\)-group is properly contained in
its normalizer, so \(H=C_G(s)\). It therefore suffices to show that
\(\alpha(s)\) is conjugate to \(s\).

The automorphism \(\alpha\) preserves each \(Pd^i\) and fixes each
element of \(Z\), leaving the case \(s=(x,z)\) with \(x\ne0\).
By \eqref{eq:commutator}, conjugation by \((xt,0)\) changes the second
coordinate of \(s\) by \(x^6(t^4+t^2)\).
The image of \(t\mapsto t^4+t^2\) is the kernel of the absolute trace
\(\Tr_{F/\mathbb F_2}\). Indeed, squaring is bijective on \(F\), and the
map \(u\mapsto u^2+u\) has kernel \(\mathbb F_2\) and image in the
trace kernel. Both the image and the trace kernel have order eight.
Since \(x^{-5}\in\mathbb F_4\), its trace from \(F\) to \(\mathbb F_4\)
is zero, so its absolute trace is also zero. Choose \(t\in F\) with
\(t^4+t^2=x^{-5}\). Conjugation by \((xt,0)\) then sends \(s\) to
\((x,z+x)=\alpha(s)\).

Now assume that \(s=1\). Then \(H=N_G(H)\). If \(H<P\), the normalizer
condition in \(P\) gives \(H<N_P(H)\). If \(H=P\), then \(N_G(H)=G\).
Thus \(H\not\le P\), so \(H\) contains a subgroup of order five.
By Sylow conjugacy, we may assume that \(D\le H\). It follows that
\(H=(H\cap P)\rtimes D\).

The action of \(D\) is irreducible on both \(Z\) and \(P/Z\). In fact,
multiplication by \(a\) on \(F\) has irreducible minimal polynomial of
degree four, since \(a\) has order five and belongs to no proper
subfield of \(F\). Formula \eqref{eq:commutator} gives \(1
e P'\le Z\), since
\([(1,0),(y,0)]
e 1\) for \(y\in F\setminus\mathbb F_2\).
Since \(P'\) is \(D\)-invariant and \(Z\) is irreducible, we have \(P'=Z\).

Let \(U=H\cap P\). If \(U\le Z\), irreducibility gives \(U=1\) or \(U=Z\).
Otherwise the image of \(U\) in \(P/Z\) is all of \(P/Z\).
In that case \(P=UZ\), and centrality of \(Z\) gives \(U'=P'=Z\).
Thus \(Z\le U\) and \(U=P\). We conclude that \(H\) is conjugate to
\(D\), \(Z\rtimes D\), or \(G\). Each of these subgroups is preserved
by \(\alpha\). Hence \((\alpha(H),1)\) is conjugate to \((H,1)\).

Proposition~\ref{prop:obstruction} now applies.
\end{proof}

\section*{Use of generative AI}

Under the author's direction and supervision, OpenAI's ChatGPT assisted
with aspects of developing the counterexample and proof, checking
calculations, and preparing the manuscript.

\end{document}